\documentclass[10pt]{amsart}
\usepackage{amssymb}
\usepackage{MnSymbol}

\usepackage{amscd}
\usepackage{amsthm}
\usepackage[dvips]{graphicx}
\usepackage{color}
\usepackage[all]{xy}
\usepackage{verbatim}

\newtheorem{Lemma1}{{Lemma}}[section]
\newtheorem{Theo1}[Lemma1]{{Theorem}}
\newtheorem*{Theo2}{{Theorem}}
\newtheorem{Def1}[Lemma1]{{Definition}}
\newtheorem{Prop1}[Lemma1]{{Proposition}}
\newtheorem{Claim1}[Lemma1]{{Claim}}
\newtheorem{Rem1}[Lemma1]{{Remark}}
\newtheorem{Cor1}[Lemma1]{{Corollary}}
\newtheorem{Ex1}[Lemma1]{{Example}}
\newtheorem{Not1}[Lemma1]{{Notation}}

\newenvironment{Lemma}{\begin{Lemma1}}{\end{Lemma1}}
\newenvironment{Def}{\begin{Def1}\rm}{\end{Def1}}
\newenvironment{Prop}{\begin{Prop1}}{\end{Prop1}}
\newenvironment{Rem}{\begin{Rem1}\rm}{\end{Rem1}}
\newenvironment{Theorem}{\begin{Theo1}}{\end{Theo1}}

\newenvironment{Cor}{\begin{Cor1}}{\end{Cor1}}

\title{Clifford's theorem for orbit categories}

\author{Alexander Zimmermann; orcid: 0000-0001-8168-4624}
\address{\newline
Universit\'e de Picardie,
\newline D\'epartement de Math\'ematiques et LAMFA (UMR 7352 du CNRS),
\newline 33 rue St Leu,
\newline F-80039 Amiens Cedex 1,
\newline France}
\email{alexander.zimmermann@u-picardie.fr}

\date{December 4, 2021; revised January 16, 2023}

\newcommand{\uar}{\uparrow}
\newcommand{\dar}{\downarrow}
\newcommand{\lra}{\longrightarrow}

\newcommand{\ra}{\rightarrow}
\newcommand{\sdp}{\rtimes} 
\newcommand{\epi}{\lra \kern-.8em\ra}

\newcommand{\dickebox}{{\vrule height5pt width5pt depth0pt}}

 \normalbaselines

\begin{document}

\begin{abstract}
Clifford theory relates the representation theory of finite groups
to those of a fixed normal subgroup by means of induction and
restriction, which is an adjoint pair of functors. We generalize
this result to the situation of a Krull-Schmidt category on which a finite
group acts as automorphisms. This then provides the orbit category
introduced by Cibils and Marcos, and studied intensively by Keller
in the context of cluster algebras, and by Asashiba in the context
of Galois covering functors. We formulate and prove Clifford's theorem
for Krull-Schmidt orbit categories with respect to a finite group $\Gamma$ of
automorphisms, clarifying this way how the image of an indecomposable object
in the original category
decomposes in the orbit category. The pair of adjoint functors
appears as the Kleisli category
of the naturally appearing monad given by $\Gamma$.
\end{abstract}

\maketitle

\section{Introduction}

Clifford theory for finite groups links the representation theory of a normal subgroup
$N$ of $G$ to the representation theory of $G$. It is known that large parts of
this classical theory does not depend on the coefficient domain $R$. The crucial
part is the notion of the inertia group $I_G(M)$ of an indecomposable $RN$-module $M$,
which is defined as the subgroup of $G$ formed by those elements $g\in G$,
such that the twisted $RN$-module $^gM$ is isomorphic to $M$ as an $RN$-module.
The most elementary part of the theory shows that
then for any indecomposable direct factor $M_0$ of the induced $RI_G(M)$-module
$M\uar_N^{I_G(M)}$ we get that $M_0\uar_{I_G(M)}^G$ is indecomposable.
A Krull-Schmidt situation really is natural, and indeed necessary, for this
statement.

Orbit categories arise in representation theory at three places at least. Cluster categories
are constructed using orbit categories of triangulated categories, as it was made precise by
Keller~\cite{triangulatedorbitcategories}. Similarly, Peng and Xiao~\cite{PengXiao1,PengXiao2}
used quotient categories for their construction of quantum groups as Hall algebras of
derived categories of hereditary algebras. Further, Riedtmann~\cite{Riedtmann},
Cibils, Solotar and Redondo~\cite{CibilsSolotar,CibilsSolotarRedondo1,CibilsSolotarRedondo2}
study Gabriel's Galois covering technique in an abstract fashion. More
systematically Cibils and Marcos~\cite{CibilsMarcos}, and later
Asashiba~\cite{AsashibaGalois1} use orbit categories to explain and
actually define clearly these categories. Further, we mention that the setting appears in
the context of braided tensor categories and fusion categories (cf~\cite{Drinfeld} and
\cite{EtingofGelaki}).

The present paper uses arguments from Clifford theory of finite group representation
to answer the question what is the decomposition of the image of indecomposable
objects in the orbit category. First, orbit categories tend to be non idempotent complete.
Hence one needs to consider the Karoubi envelope of the orbit categories.
Then, we use techniques from category theory, namely the Kleisli adjunction of a
monad, to give an
analogue of Clifford's theory for orbit categories.

In a sense the present result can also be seen as a continuation of our previous
results~\cite{Greentriangulated,GreenSEAMS}
on a categorical framework for Green correspondence for adjoint functors.
Clifford theory and Green correspondence are two of the main building blocks of
modular representation theory of finite groups. Giving a framework in
a category framework, opening hence the possible applications is a
desirable task. Auslander-Kleiner
\cite{AuslanderKleiner} showed that classical Green correspondence can be
formulated and proved as a property of a pair of adjoint functors between
additive categories. This was then put in the framework of triangulated categories
in \cite{Greentriangulated,GreenSEAMS}. Alternative approaches were given by
Benson, Carlson, Grime, Peng, Wheeler, Wang and Zhang in a series of joint papers
\cite{Carlson-Peng-Wheeler,Grime,Benson-Wheeler,Carlson,Wang-Zhang}. Another different
approach is given by Balmer and del-Ambrogio \cite{BalmerDellAmgrogio} in the context of
tensor triangulated categories and Mackey $2$-functors.

We start with a Krull-Schmidt category $\mathcal H$ on which a finite group $\Gamma$ acts.
The orbit category ${\mathcal G}:={\mathcal H}[\Gamma]$ is then actually the
Kleisli construction of the corresponding
monad and provides a pair of adjoint functors $(S,T)$
$$\xymatrix{{\mathcal H}\ar@/^/[r]^S&{\mathcal G}\ar@/^/[l]^T}$$
such that $TS=\bigoplus_{i\in I}E_i$ is the direct sum of
automorphisms $E_i$ of $\mathcal H$, which forms the group $\Gamma$.
We replace the inertia group in the group algebra situation
by the orbit category ${\mathcal H}[\Gamma_M]$
with respect to the subgroup $\Gamma_M$ of those elements of
$\Gamma$ which fix the isomorphism class of a given object $M$. We first show that
$\Gamma$ lifts to a group of automorphisms of $\mathcal G$ and of
${\mathcal H}[\Gamma_M]$.
We obtain that
the lift $\widehat \Gamma$ of $\Gamma$ as group of automorphisms of $\mathcal G$
fixes isomorphism classes of objects.

These orbit categories  ${\mathcal H}[\Gamma_M]$ and $\mathcal G$  do not
have split
idempotents in general, but replacing them with their Karoubi envelope, we can
show a precise analogue of Clifford's theorem. This is our main
result Theorem~\ref{CliffordforMonads}.
We note that our main Theorem~\ref{CliffordforMonads} starts with an
$R$-linear Krull-Schmidt category and an action of a finite group $\Gamma$
on $\mathcal H$. Then Theorem~\ref{CliffordforMonads} determines how
indecomposable objects in $\mathcal H$ behave in the orbit
category ${\mathcal H}[\Gamma]$ with respect to indecomposability.

The paper is organised as follows. In Section~\ref{monadsection} we recall
the necessary notations and background on monads, the Kleisli category,
the Eilenberg-Moore adjunctions and their properties as far as we need them.
The short Section~\ref{Groupsituationsection} recalls the group algebra
situation.
In Section~\ref{inertiacasesection} we define categories which model
a normal subgroup in our setting, and we prove our first main result,
Theorem~\ref{Resultincaseoftrivialinertiagroup}, which is the Clifford theorem
for monads in the case when analogue of
the inertia subgroup does not exceed the normal subgroup.
This situation does not need all the hypotheses, and is therefore formulated
in a more general setting. Section~\ref{generalcasesection}
then defines an analogue of the inertia group in our general situation,
studies its properties and shows the main result Theorem~\ref{CliffordforMonads}.
We use in particular properties on orbit categories, Karoubi envelopes and
some general statements on adjoint functors. All these tools are recalled
in this section. Finally in Section~\ref{Exampleapplicationsection} we present
examples from Galois modules and from fusion categories.

\section{Monads revisited}
\label{monadsection}

Recall that if $(S,T)$ is an adjoint pair, then the endofunctor
$TS$ of $\mathcal H$ together with the unit and the counit of the
adjunction give a monad. There is extensive literature on monads.

\begin{Def} \cite{Maclane}
A monad $(A,\mu,\eta)$ on a category $\mathcal C$ is an endofunctor $A$ of $\mathcal C$
with a natural transformation $$\mu:A^2\lra A$$ such that
$$\mu\circ (A\mu)=\mu\circ(\mu A):A^3\lra A$$
and $$\eta:\textup{id}_{\mathcal C}\lra A$$ such that $$\mu\circ(A\eta)=\mu\circ(\eta A).$$

If ${\bf A}:=(A,\mu,\eta)$ is a monad, an ${\bf A}$-module in $\mathcal C$ is a pair
$(X,\rho)$ where $X$ is an object in $\mathcal C$ and $\rho:A(X)\lra X$ is a morphism
in $\mathcal C$ such that $$\rho\circ(A\rho)=\rho\circ\mu_X:A^2(X)\lra X$$
and $$\rho\circ\eta_X=\textup{id}_X:X\lra X.$$

If ${\bf A}:=(A,\mu,\eta)$ is a monad and $(X,\rho)$ and $(X',\rho')$ are $\mathbf A$-modules, then
a morphism of $\mathbf A$-modules $$f:(X,\rho)\lra (X',\rho')$$
is an $f\in{\mathcal C}(X,X')$ such that $$\rho'\circ A(f)=f\circ\rho:A(X)\lra X'.$$

${\mathbf A}-\textup{Mod}_{\mathcal C}$ is the category (!) of ${\mathbf A}$-modules on $\mathcal C$.
\end{Def}

Every adjunction induces a monad. In general many different adjunctions induce
the same monad and all these form a
category. There are two particular and extremal such adjunctions realising a given monad.
The Eilenberg-Moore adjunction is a
terminal object in this category of adjunctions realising a give monad
and the Kleisli category is an initial object
in this category  of adjunctions realising a give monad.

\subsection{The Eilenberg-Moore adjunction}

\label{EilenbergMooreSection}

The Eilenberg-Moore \cite{EilenbergMoore} adjunction is an adjoint pair
$(F_{\mathbf A},U_{\mathbf A})$ where
$F_{\mathbf A}:{\mathcal C}\lra {\mathbf A}-\textup{Mod}_{\mathcal C}$
is defined as $$F_{\mathbf A}(Y):=(A(Y),\mu_Y)\textup{ and }F_{\mathbf A}(f)=A(f)$$
and $$U_{\mathbf A}(X,\rho)=X\textup{ and }U_{\mathbf A}(f)=f.$$
The objects in the image of $F_{\mathbf A}$ are called free $\mathbf A$-modules, and the full
subcategory of ${\mathbf A}-\textup{Mod}_{\mathcal C}$ generated by free modules is denoted
${\mathbf A}-\textup{Free}_{\mathcal C}$. Restricting to the image of $F_{\mathbf A}$
one obtains the so-called Kleisli adjunction
$(\widehat F_{\mathbf A},\widehat U_{\mathbf A})$ where
$\widehat F_{\mathbf A}:{\mathcal C}\lra {\mathbf A}-\textup{Free}_{\mathcal C}$.

If $(S,T)$ is an adjoint pair with counit $\epsilon:ST\lra \textup{id}_{\mathcal D}$,
and where $S:{\mathcal C}\lra{\mathcal D}$, then $TS$ is a monad with $\mu=T\epsilon S$
and $\eta:{\textup id}_{\mathcal C}\lra TS$ is the unit.

Further, by \cite[Chapter VI]{Maclane} if $(S,T)$ is an adjoint pair with
$S:{\mathcal C}\lra{\mathcal D}$ we get functors
$$K:{\mathbf A}-\textup{Free}_{\mathcal C}\lra{\mathcal D}\textup{ and }
E:{\mathcal D}\lra {\mathbf A}-\textup{Mod}_{\mathcal C}$$
such that
$$\widehat U_{\mathbf A}=T\circ K ,\;
S=K\circ \widehat F_{\mathbf A},\; F_{\mathbf A}=E\circ S, \;
T= U_{\mathbf A}\circ E.$$
and in particular the diagram $(\dagger)$
$$
\xymatrix{
&{\mathcal C}\ar[d]^S\ar[dl]_{\widehat F_{\mathbf A}}\ar[dr]^{ F_{\mathbf A}}\\
{\mathbf A}-\textup{Free}_{\mathcal C}\ar[r]_-K&{\mathcal D}\ar[r]_-E&{\mathbf A}-\textup{Mod}_{\mathcal C}
}
$$
is commutative, and such that $E\circ K$ is a fully faithful embedding.
Here, $K$ is defined by
\begin{eqnarray*}
{\mathbf A}-\textup{Free}_{\mathcal C}&\lra&{\mathcal D}\\
F_{\mathbf A}(Y)&\mapsto&S(Y)
\end{eqnarray*}
and by
\begin{eqnarray*}
{\mathbf A}-\textup{Free}_{\mathcal C}(F_{\mathbf A}(Y),F_{\mathbf A}(Y'))
&\simeq&{\mathcal C}(Y,AY')\\
&\simeq&{\mathcal D}(SY,SY')\\
&\simeq&{\mathcal D}(KF_{\mathbf A}Y,KF_{\mathbf A}Y')
\end{eqnarray*}
on morphisms. Hence $K$ is always fully faithful.
Further, $E(Z)=(T(Z),T(\epsilon_Z))$ on objects and $E(f)=T(f)$ on morphisms.
An adjunction $(S,T)$ is called {\em monadic} if $E$ is an equivalence.
A monad $(A,\mu,\eta)$ on a category $\mathcal C$ is called {\em separable}
if there is a natural transformation $\sigma:A\lra A^2$ such that
$\mu\circ\sigma=\textup{id}_A$ and
$$(A\mu)\circ(\sigma A)=\sigma\circ\mu=(\mu A)\circ(A\sigma):A^2\lra A^2.$$

\begin{Rem}
In the following Proposition~\ref{monadsadmittingsection}
Balmer uses the concept of an equivalence up to direct summands.
A functor $F:{\mathcal C}\lra {\mathcal D}$ is called to be an equivalence up to direct summands if
the induced functor $\widehat{F}:Kar({\mathcal C})\lra Kar({\mathcal D})$ is an equivalence.
For more details on the Karoubi idempotent completion $Kar$ and the functor $\widehat F$ see
Remark~\ref{Karoubienvelope} below.
\end{Rem}

\begin{Prop}\cite[Lemma 2.10]{BalmerStacks} \label{monadsadmittingsection}
Let $(S,T)$ be a pair of adjoint
functors, $S:{\mathcal C}\lra{\mathcal D}$ such that the counit
$\epsilon:ST\lra \textup{id}_{\mathcal D}$ has a section
$\xi:\textup{id}_{\mathcal D}\lra ST$, i.e.
$\xi\circ\epsilon=\textup{id}$. Then
\begin{enumerate}
\item The monad $A=TS$ is separable
\item The functors $K$ and $E$ are equivalences up to direct summands
\item If $\mathcal C$ and $\mathcal D$ are idempotent complete, then
$K$ and $E$ are equivalences.
\end{enumerate}
\end{Prop}

\subsection{The Kleisli category}

Kleisli gave another direct construction realising the monad by an adjunction, the Kleisli category,
which is actually isomorphic to ${\mathbf A}-\textup{Free}_{\mathcal C}$.

Given a monad $(A,\mu,\eta)$ on $\mathcal H$,
then define the Kleisli category ${\mathfrak K}_{A}$ by the following construction:

The objects of ${\mathfrak K}_{A}$ and of $\mathcal H$ coincide.
Let $X,Y$ be two objects. Then
$$
{\mathfrak K}_{A}(X,Y):=\{f\in{\mathcal H}(AX,AY)\;|\;\xymatrix{
A^2X\ar[r]^{Af}\ar[d]_{\mu_{AX}}&A^2Y\ar[d]^{\mu_{AY}}\\ AX\ar[r]_{f}&AY
}\textup{ is commutative}\}
$$
Composition is given by composition of maps in $\mathcal H$. This is well-defined
since $A$ is a functor. Then there are functors
$$T_{{\mathfrak K}_{A}}:{\mathfrak K}_{A}\lra {\mathcal H}$$
and $$S_{{\mathfrak K}_{A}}:{\mathcal H}\lra {\mathfrak K}_{A}$$  given by
the following:

$T_{{\mathfrak K}_{A}}(X)=AX$ for any object $X$ and $T_{{\mathfrak K}_{A}}(f)=f$
for morphisms.

$S_{{\mathfrak K}_{A}}(X)=X$ for any object $X$ and $S_{{\mathfrak K}_{A}}(f)=Af$
for morphisms.

Again, $(S_{{\mathfrak K}_{A}},T_{{\mathfrak K}_{A}})$
is an adjoint pair inducing the monad $A$ (cf e.g. \cite[Chapter 2]{Pareigis}. )

\medskip

The following lemma seems to be well-known to the specialist
(cf e.g. the introduction into Teleiko \cite{Teleiko})
but for the convenience of the reader we include the result and the short proof.

\begin{Lemma} \label{homsforkleisli}
Let $\mathcal H$ be a category and let $(A,\mu,\eta)$ be a monad on $\mathcal H$.
Then ${\mathfrak K}_{A}(X,Y)={\mathcal H}(X,AY)$.
\end{Lemma}

Proof. We have the unit $\eta:\textup{id}_{\mathcal H}\lra A$ and get a map
\begin{eqnarray*}
{\mathfrak K}_A(X,Y)&\stackrel{\phi}\lra&{\mathcal H}(X,AY)\\
f&\mapsto&\eta_X\circ f
\end{eqnarray*}
This way the diagram
$$
\xymatrix{
AX\ar[r]^f&AY\\
X\ar[u]_{\eta_X}\ar[ru]_{\phi(f)}
}
$$
is commutative. We have a map in the opposite direction.
\begin{eqnarray*}
{\mathcal H}(X,AY)&\stackrel{\psi}\lra&{\mathcal H}(AX,AY)\\
g&\mapsto&\mu_Y\circ A(g)
\end{eqnarray*}
$$
\xymatrix{&A^2Y\ar[d]^{\mu_Y}\\
AX\ar[ru]^{Ag}\ar[r]_{\psi(g)}&AY
}
$$
We compose $\phi\circ\psi$ as in the diagram
$$
\xymatrix{&A^2Y\ar[d]^{\mu_Y}\\
AX\ar[ru]^{Ag}\ar[r]_{\psi(g)}&AY\\
X\ar[u]_{\eta_X}\ar[ru]_{\phi(\psi(g))}
}
$$
But since $\mu\circ A\eta=\mu\circ\eta A=\textup{id}$,
$$\mu_Y\circ A(g)\circ\eta_X=\mu_Y\circ \eta_{A(Y)}\circ g=g.$$
Hence $\phi\circ\psi=\textup{id}.$

In order to show that $\psi\circ\phi$ is the identity
we consider $f\in {\mathfrak K}_A(X,Y)$.
$$
\xymatrix{
AX\ar[r]^f&AY\\
X\ar[u]^{\eta_X}\ar[ur]_{\phi(f)}
}
$$
gives the commutative diagram
$$
\xymatrix{
A^2X\ar[rr]^{Af}&&A^2Y\ar[d]^{\mu_Y}\\
AX\ar[u]^{A\eta_X}\ar[urr]_{A\phi(f)}\ar[rr]_{\psi(\phi(f))}&&AY
}
$$
But since $f\in {\mathfrak K}_A(X,Y)$, the diagram
$$
\xymatrix{
A^2X\ar[rr]^{Af}\ar[d]_{\mu_X}&&A^2Y\ar[d]^{\mu_Y}\\
AX\ar[rr]_{f}&&AY
}
$$
is commutative. Hence
$$\psi\phi(f)=\mu_Y\circ Af\circ A\eta_X=f\circ\mu_X\circ A\eta_X=f$$
since $\mu\circ A\eta=\textup{id}$. \dickebox

\medskip

Note that in order to compute the composition of two maps, the original
definition of morphisms in the Kleisli construction is more appropriate.

\subsection{Universal  property of the two Eilenberg-Moore and the Kleisli category}

\label{KleisliEMsection}

Let $(A,\mu,\eta)$ be a monad on $\mathcal H$. Let moreover
$$
\xymatrix{{\mathcal H}\ar@/^/[r]^S&{\mathcal G}\ar@/^/[l]^T}
$$
be functors such that $(S,T)$ is an adjoint pair inducing the monad $(A,\mu,\eta)$ with $A=TS$.
Then there are unique functors $K:{\mathfrak K}_{A}\lra {\mathcal G}$
and $L:{\mathcal G}\lra A-Mod_{\mathcal H}$ making the diagram
$$
\xymatrix{
&&A-Mod_{\mathcal H}\ar@/_/[ddll]|{T_{EM}}\\ \\
{\mathcal H}\ar@/_/[uurr]|{S_{EM}}\ar@/^/[rr]|S\ar@/^/[ddrr]|{S_{{\mathfrak K}_{A}}}&&{\mathcal G}\ar@/^/[ll]|{T}\ar[uu]|L\\ \\
&&{\mathfrak K}_{A}\ar@/^/[uull]|{T_{{\mathfrak K}_{A}}}\ar[uu]|{K}
}
$$
commutative, in the sense that
$$S=KS_{{\mathfrak K}_{A}}; T_{{\mathfrak K}_{A}}=TK$$
$$S_{EM}=LS; T=T_{EM}L.$$

Indeed, on objects we have:
$K(X):=S(X)$ for all objects $X$ of ${\mathfrak K}_{A}$, which are the same
objects as those in $\mathcal H$.
As for morphisms we obtain
$$
\xymatrix{
{\mathfrak K}_{\mathcal H}(X,Y)\ar@{^{(}->}[r]
&{\mathcal H}(AX,AY)\ar[r]^-{\simeq}
&{\mathcal G}(STSX,SY)\ar[r]& {\mathcal G}(SX,SY)\\
f\ar@{|->}[rrr]&&&L(f)
}
$$
where the first isomorphism is the adjointness property, and the latter map
is given by the unit $\eta:\textup{id}_{\mathcal G}\lra TS$, which induces a
natural transformation $S\eta:S\lra STS$, which is actually split,
by the fundamental property of adjoint functors.
For any morphism $f\in {{\mathfrak K}_{\mathcal H}}(X,Y)$ put
$K(f)$ the image under this map.

As for $L$ we put $L(X):=(T(X),\epsilon_Y)$ on objects,
where $\epsilon:ST\lra \textup{id}_{\mathcal H}$ is the unit of the adjointness,
and $L(f):=T(f)$ on morphisms.

Note that the above property, together with the diagram $(\dagger)$
in Section~\ref{EilenbergMooreSection}, show that the Kleisli category and
${\mathbf A}-\textup{Free}_{\mathcal C}$ are equivalent.
For more details see e.g. \cite[Section 2.3 Satz 1]{Pareigis}.

\subsection{A further property of adjoint functors}

\label{unitcounitsplits}
We shall need a further property of adjoint functors. These statements
are well-known and
can be extracted implicitly from the proof of \cite[Chapter 1, Proposition 1.3]{GabrielZisman}.
For the convenience of the reader we provide a short proof.

\begin{Lemma}\label{adjointpairsfullandfaithful}
Let ${\mathcal H}$ and ${\mathcal G}$ be categories and let
$$
\xymatrix{{\mathcal H}\ar@/^/[r]^S&{\mathcal G}\ar@/^/[l]^T}
$$
be functors such that $(S,T)$ is an adjoint pair.
Then
\begin{enumerate}
\item $S$ is faithful if and only if the unit $\textup{id}\lra TS$ is pointwise a monomorphism.
\item $T$ is faithful if and only if the counit $ST\lra \textup{id}$ is pointwise an epimorphism.
\item $S$ is full if and only if the unit $\textup{id}\lra TS$ is pointwise a split epimorphism.
\item $T$ is full if and only if the counit $ST\lra \textup{id}$ is pointwise a split monomorphism.
\item $S$ is fully faithful if and only if the unit $\textup{id}\lra TS$ is an isomorphism.
\item $T$ is fully faithful if and only if the counit $ST\lra \textup{id}$ is an isomorphism.
\end{enumerate}
\end{Lemma}

Proof. Items $(2n)$ and $(2n-1)$ are dual to each other for all $n\in\{1,2,3\}$.
We need to show one of the corresponding statements.

Item 2):
$T$ is faithful if and only if ${\mathcal G}(X,Y)\stackrel T\lra {\mathcal H}(TX,TY)$ is injective.
Now ${\mathcal H}(TX,TY)\simeq {\mathcal G}(STX,Y)$ and any map in the latter factors through
the counit $\epsilon_X:STX\lra X$ (cf \cite[IV.1 Theorem 1]{Maclane}). Now,
$\epsilon_X$ is an epimorphism if and only if
${\mathcal G}(\epsilon_X,Y)$ is injective and this shows the statement.

Item 4):
$T$ is full if and only if
${\mathcal G}(X,Y)\stackrel T\lra {\mathcal H}(TX,TY)$ is surjective.
Again ${\mathcal H}(TX,TY)\simeq {\mathcal G}(STX,Y)$ and any map in the
latter factors through
the counit $\epsilon_X:STX\lra X$ (cf \cite[IV.1 Theorem 1]{Maclane}).
Hence $T$ is full if and only if
${\mathcal G}(X,Y)\stackrel{{\mathcal G}(\epsilon_X,Y)}{\lra} {\mathcal G}(STX,Y)$
is surjective
for all $X,Y$. In particular, for $Y=STX$ we obtain that
this is equivalent to $\epsilon_X$ being split monomorphism.

Item 6): An epimorphism which is in addition a split monomorphism is an isomorphism.
\dickebox

\section{Recall Clifford's theorem from group algebras}

\label{Groupsituationsection}

We briefly recall the situation for group rings.
Let $G$ be a group and let $k$ be a commutative ring.
For a subgroup $H$ of $G$ and a $kH$-module $M$ we denote by $M\uar_H^G$ the
$kG$-module $kG\otimes_{kH}M$, where $kG$ acts by multiplication on the left factor.
Further, for a $kG$-module $X$ we denote by $X\dar^G_H$ the
$kH$-module obtained by restricting the action of $G$ to an action of the subgroup $H$.
Both $\uar_H^G$ and $\dar^G_H$ are functors between the respective module categories, and
actually the are biadjoint one to the other in case $G$ is a finite group.

If $N\unlhd G$, then  Mackey's formula shows
$$M\uar_N^G\dar^G_H=\bigoplus_{gN\in G/N }\ ^gM$$
for each $M$, and actually we have an isomorphism of functors
$$ -\uar_N^G\dar^G_N=\bigoplus_{gN\in G/N }\ ^g-.$$

\begin{Rem}\label{directsumuniversal}
We should note that $M\uar_H^G=kG\otimes_{kH}M$ and $N\dar^H_G=Hom_{kH}(kG,N)$ both only depend
on the group algebras, rather than the groups. Therefore, the decomposition
is not canonical, and since there may be different groups with the same group algebra,
the decomposition depends on the specific group one takes.
We mention \cite{Margolis} for a striking example in characteristic $3$ and \cite{Hertweck} for
an example over the integers.
\end{Rem}

If $N$ is a normal subgroup of $G$ and $M$ an indecomposable $kN$-module. Then
the inertia group $I_G(M)$ of $M$ is the set of $g\in G$ such that $^gM\simeq M$
as $kN$-module $$I_G(M):=\{g\in G\;|\;^gM\simeq M\textup{ as $kN$-modules}\}.$$

\begin{Theorem} (Clifford \cite{Clifford})
\label{Cliffordtheoremclassic}
Let $k$ be a field and let $G$ be a group with normal subgroup $N\unlhd G$ of
finite index.
Let $M$ be an indecomposable $kN$-module and let
$M_0$ be an indecomposable direct factor of $M\uar_N^{I_G(M)}$. Then
$M_0\uar_{I_G(M)}^G$ is an indecomposable $kG$-module.
\end{Theorem}

In the rest of the paper we show that this theorem is actually a result on orbit categories.

\section{Categories modelling normal subgroups}

\label{inertiacasesection}

Motivated by Mackey's theorem for group rings in the situation of normal subgroups
we give

\begin{Def}
Let ${\mathcal H}$ and ${\mathcal G}$ be categories and let
$$
\xymatrix{{\mathcal H}\ar@/^/[r]^S&{\mathcal G}\ar@/^/[l]^T}
$$
be functors such that $(S,T)$ is an adjoint pair. Then $S$ gives
a {\em situation of normal subgroup categories} if
$$A=TS=\bigoplus_{i\in I}E_i$$
for self-equivalences $E_i$; $i\in I$; of $\mathcal H$.
\end{Def}

Note that in the main result of Auslander and Kleiner \cite{AuslanderKleiner}
as well as in its triangulated category version \cite{Greentriangulated}
we started from the situation
that $(S,T)$ is an adjoint pair and $TS={\textup{id}}_{\mathcal H}\oplus U$
for some endofunctor $U$ of $\mathcal H$, and such that the unit
$\textup{id}_{\mathcal H}\lra TS$ is the left inverse to the projection
$TS\lra \textup{id}_{\mathcal H}$. Hence, if $\mathcal H$ is Krull-Schmidt,
we may assume that for some index $i_0\in I$ we have $E_{i_0}=\textup{id}_{\mathcal H}$.
Note further that by Remark~\ref{directsumuniversal} it is not reasonable to assume the
decomposition to be canonical in any sense. We only ask for the existence.

A particularly simple situation occurs in the classical case
of Clifford's theorem if $N=I_G(M)$.
In our categorical situation this is a favorable situation as well.

\begin{Prop}
Let $S:{\mathcal H}\lra{\mathcal G}$ be an additive functor between
abelian Krull-Schmidt categories admitting a right
adjoint $T$.
Assume
further that $S$ is a normal subgroup category situation. We have, by definition
$TS=\bigoplus_{i\in I} E_i$ for self-equivalences $E_i$ of $\mathcal H$,
and for some $i_0$, $E_{i_0}=\textup{id}$.
If $E_iM\not\simeq M$ for any $i\in I\setminus\{i_0\}$, then
$M\textup{ simple }$ implies that $End_{\mathcal G}(SM)$ is a skew field and in particular
$SM\textup{ is indecomposable }$.
\end{Prop}

Proof.
We compute
\begin{eqnarray*}
End_{\mathcal G}(SM)&\simeq& {\mathcal H}(M,TSM)\\
&\simeq&
{\mathcal H}(M,E_{i_0}M)\oplus\bigoplus_{i\in I\setminus\{i_0\}}{\mathcal H}(M,E_{i}M)
\end{eqnarray*}
Since $E_i$ are all self-equivalences, $M$ is simple if and only if $E_iM$ is simple.
Hence, ${\mathcal H}(M,E_{i}M)=0$ if $i\in I\setminus\{i_0\}$.
Further, by functoriality of the adjointness, the resulting isomorphism
$End_{\mathcal G}(SM)\simeq End_{\mathcal H}(M)$ is an isomorphism of rings.
Therefore $End_{\mathcal G}(SM)$ is a skew field and therefore $SM$ is indecomposable.
\dickebox

\medskip

Suppose that $E_iA\stackrel{\lambda_i}{\lra}AE_i$ are natural equivalences.
Then for all $f\in {\mathcal H}(X,Y)$ we have a commutative diagram
$$ (\dagger)
\xymatrix{
E_iAX\ar[r]^-{E_iAf}\ar[d]_{(\lambda_i)_X}&E_iAY\ar[d]^{(\lambda_i)_Y}\\
AE_iX\ar[r]^-{AE_if}&AE_iY}
$$
where the vertical maps are isomorphisms. As a consequence, if $E_i$ is an
automorphism, then
$$AE_i^{-1}\stackrel{E_i^{-1}(\lambda_i)E_i^{-1}}{\lra}E_i^{-1}A$$
is a natural equivalence in the sense that
$$ (\dagger\dagger)
\xymatrix{
AE_i^{-1}X\ar[r]^-{AE_i^{-1}f}\ar[d]_{E_i^{-1}(\lambda_i)_{E_i^{-1}X}}&
AE_i^{-1}Y\ar[d]^{{E_i^{-1}(\lambda_i)_{E_i^{-1}Y}}}\\
E_i^{-1}X\ar[r]^-{E_i^{-1}Af}&E_i^{-1}AY}
$$
is commutative.

\begin{Prop}\label{EiliftstoG}
Let $S:{\mathcal H}\lra{\mathcal G}$ be an additive functor between
Krull-Schmidt categories admitting a right adjoint $T$. Assume
further that $S$ is a normal subgroup category situation. Then
$TS=\bigoplus_{i\in I} E_i$ for self-equivalences $E_i$ of $\mathcal H$.
Further, assume that $E_{i_0}=:E_0=\textup{id}_{\mathcal H}$.
Denote by $A=TS$ the endofunctor of the monad ${\mathbf A}$
induced by this adjoint pair.
Let $J\subseteq I$ and suppose that for each $j\in J$
there are natural equivalences $E_jA\stackrel{\lambda_j}{\lra}AE_j$.
Then for each $j\in J$
there is a self-equivalence of the Kleisli category
$\widehat E_j$ of ${\frak K}_{A}$ such that
$E_j\circ T=T\circ \widehat E_j$ and $\widehat E_j\circ S=S\circ E_j$. If $E_j$
is an automorphism, then
$\widehat E_j$ is an automorphism. Analogous statements hold
for the Eilenberg-Moore category
$A-\textup{Mod}_{\mathcal H}$.
\end{Prop}

Proof. Let $j\in J$.
We shall prove that each $E_j$ lifts to a
self-equivalence $\widehat E_j$ of ${\frak K}_A$.

We put $\widehat E_j(X)=E_j(X)$ for each object $X$ of
${\frak K}_A$.
As for morphisms we may use Lemma~\ref{homsforkleisli}. There we obtained
${\mathcal H}(X,AY)={\mathfrak K}_{\mathbf A}(X,Y)$ and for $f\in {\mathcal H}(X,AY)$
we have $E_jf\in{\mathcal H}(E_jX,E_jAY)$, and therefore
$(\lambda_j)_Y\circ E_jf\in {\mathcal H}(E_jX,AE_jY)$. Hence we may put
$\widehat E_j(f):=(\lambda_j)_Y\circ E_jf$. If $E_j$ is an automorphism, then
the diagrams $(\dagger\dagger)$ and $(\dagger)$ show that $\widehat E_j$ is an
automorphism as well.
Indeed, $$\widehat{E_j^{-1}}(f):={E_j^{-1}(\lambda_j)_{E_j^{-1}Y}}^{-1}\circ E_j^{-1}(f)$$
gives
\begin{eqnarray*}
\widehat{E_j^{-1}}(\widehat E_j(f))&=&\widehat{E_j^{-1}}((\lambda_j)_Y\circ E_jf)\\
&=&{E_j^{-1}(\lambda_j)_{E_j^{-1}Y}}^{-1}\circ E_j^{-1}(((\lambda_j)_Y\circ E_jf)\\
&=&{E_j^{-1}(\lambda_j)_{E_j^{-1}Y}}^{-1}\circ E_j^{-1}(\lambda_j)_{E_{j}^{-1}Y}\circ f\\
&=&f
\end{eqnarray*}

On the level of ${\mathbf A}-\textup{Mod}_{\mathcal H}$ the construction is as follows.
We put $\rho_{E_jX}:=E_j(\rho_X)\circ(\lambda_j)_X^{-1}$ for every $X$
and we put $\widehat E_j(X,\rho_X):=(E_jX,\rho_{E_jX})$. If
$\rho_X$ is a splitting of the unit $\eta_X$, then $\rho_{E_jX}$ is a
splitting of $\eta_{E_jX}$. Hence, this is indeed an object of
${\mathbf A}-\textup{Mod}_{\mathcal H}$. Moreover, for
any $f\in {\mathcal H}(X,Y)$, we have that $E_j(f)\in {\mathcal H}(E_jX,E_jY)$,
and if $f$ gives rise to a morphism in ${\mathbf A}-\textup{Mod}_{\mathcal H}$,
then
$$
\xymatrix{
AX\ar[r]^-{Af}\ar[d]_{\rho_X}&AY\ar[d]^{\rho_Y}\\
X\ar[r]^-f&Y
}
$$
is commutative, and therefore
$$
\xymatrix{
E_jAX\ar[r]^-{E_jAf}\ar[d]_{E_j\rho_X}&E_jAY\ar[d]^{E_j\rho_Y}\\
E_jX\ar[r]^-{E_jf}&E_jY
}
$$
is commutative as well. Now, recall the diagram $(\dagger)$ and get a commutative diagram
$$
\xymatrix{
AE_jX\ar[r]^-{AE_jf}\ar[d]_{(\lambda_{j})_X^{-1}}&AE_jY\ar[d]^{(\lambda_{j})_Y^{-1}}\\
E_jAX\ar[r]^-{E_jAf}\ar[d]_{E_j\rho_X}&E_jAY\ar[d]^{E_j\rho_Y}\\
E_jX\ar[r]^-{E_jf}&E_jY.
}
$$
Since we were putting $\rho_{E_jX}=E_j(\rho_X)\circ(\lambda_j)_X^{-1}$ for
every $X$ we get that
then also $E_jf$ gives rise to a morphism in $A-\textup{Mod}_{\mathcal H}$.

The fact that $\widehat E_j$ is an equivalence is clear, since the construction
holds as well for $E_j^{-1}$, producing a functor $\widehat{ E^{-1}_j}$, which
equals
$\widehat E_j^{-1}$ as is verified by composing with $\widehat E_j$.

Let us verify that $S_K E_j=\widehat E_j S_K$. Indeed, $S_K(X)=X$ for any object $X$ and $S_K(f)=Af$ for any morphism $f$.
Then $S_K E_j(X)=\widehat E_j S_K(X)$ for any object.
Further, $$\widehat E_jS_Kf=\widehat E_jAf=\lambda_j^{-1}\circ E_jAf=\lambda_j^{-1}\lambda_jAE_jf=AE_jf.$$
Similarly, the facts that $T_K\widehat E_j=E_jT_K$, and that $T_{EM}\widehat E_j=E_jT_{EM}$
and
$S_{EM} E_j=\widehat E_jS_{EM}$ follow by direct inspection.
\dickebox

\begin{Rem}\label{liftingthegroup}
If in the notation and under the hypotheses of Proposition~\ref{EiliftstoG}
we have in addition $E_i\circ E_j=E_k$ for some $i,j,k\in I$, and if in addition
$$\lambda_k=\lambda_i\circ E_i(\lambda_j)$$
then we also
get $\widehat E_i\circ \widehat E_j=\widehat E_k$. This follows by the
explicit construction of $\widehat E_i$ in each of the two cases.
The equation on the isomorphisms $\lambda_*$ comes from the necessity that the
diagram
$$
\xymatrix{
E_iE_jA\ar[r]^-{E_i(\lambda_j)}&E_iAE_j\ar[r]^-{\lambda_j}&AE_iE_j\\
E_kA\ar@{=}[u]\ar[rr]^{\lambda_k}&&AE_k\ar@{=}[u]
}
$$
commutes. Note that $\lambda_*$ is a relation of vanishing of a
$1$-cocycle type. (Note
that we did not assume yet that the $\{E_i\;|\;i\in I\}$ forms a group.)
In general, for freely chosen $\lambda_i$ for $i\in I$, the
$1$-cocycle above is not vanishing. However, by the usual argument it is sufficient
that the $1$-cocycle above is a $1$-coboundary, that is there is $\alpha$ such that
$\lambda_i=\alpha^{-1}\circ E_i(\alpha)$ for all $i\in I$. Then, by an elementary computation
the vanishing of the $1$-cocycle follows.

Further, if $\{E_i\;|\;i\in I\}$ is a group, then $E_iA\lra AE_i$ is actually
the conjugation action on the summands of $A$. Since the conjugation action
is associative, i.e. $(E_iE_j)A(E_iE_j)^{-1}=E_i(E_jAE_j^{-1})E_i^{-1}$,
if $\{E_i\;|\;i\in I\}$ is a group, then $E_i\circ E_j=E_k$ implies
$\widehat E_i\circ \widehat E_j=\widehat E_k$.
\end{Rem}

\medskip

\begin{Prop}\label{widehatEiMisisotoM}\label{ehatmapstoisomorphiccopies}
Let $S:{\mathcal H}\lra{\mathcal G}$ be an additive functor between
Krull-Schmidt categories admitting a right adjoint $T$. Assume
further that $S$ is a normal subgroup category situation, i.e.
$TS=\bigoplus_{i\in I} E_i$ for self-equivalences $E_i$ of $\mathcal H$.
Further, assume that $E_{i_0}=:E_0=\textup{id}_{\mathcal H}$.
Denote by $A=TS$  the endofunctor of the monad ${\mathbf A}$
induced by this adjoint pair. Then the counit
$\epsilon:S_KT_K\lra \textup{id}_{{\mathfrak K}_A}$ of the Kleisli
category associated to ${\mathbf A}$ splits
pointwise, and we get $\widehat E_i(M)\simeq M$ for
all indecomposable objects $M$ of ${\mathfrak K}_A$ and the
lift $\widehat E_i$ of $E_i$ to
${\mathfrak K}_A$ following Proposition~\ref{EiliftstoG}.
\end{Prop}

Proof.
In order to prove this statement we first analyse the monad $(A,\mu,\epsilon)$
in the normal subgroup situation, in particular $\mu:A^2\lra A$.
Recall the construction of $\mu$.
Let $\epsilon:ST\lra \textup{id}_{\mathcal G}$ be the counit of the adjunction.
Then $$\mu=T\epsilon S:TSTS\lra T \circ\textup{id}_{\mathcal G} \circ S.$$
We study this map in particular for the Kleisli category.
Then, for the category ${\mathfrak K}_A$ and the adjoint pair $(S_K,T_K)$ giving the
monad $A=TS$, we have
$$\begin{array}{rclcrcl}
S_K:{\mathcal H}&\lra&{\mathfrak K}_A&&T_K:{\mathfrak K}_A&\lra&{\mathcal H}\\
X&\mapsto&X&&X &\mapsto&AX\\
f&\mapsto&Af&&f&\mapsto&f
\end{array}$$
where $f$ denotes a morphism and $X$ an object of $\mathcal H$.
Hence, $S_KT_K(X)=AX$ and $S_KT_K(f)=Af$ for all objects $X$ and all morphisms $f$.
The choice $\epsilon:S_KT_K\lra \textup{id}_{{\mathfrak K}_A}$ given by $\eta$ the projection
of $A=1\oplus U$ onto $1$ satisfies the necessary properties. Moreover, by definition
the counit defined in this way splits pointwise.

The fact that $M\simeq \widehat E_iM$
in ${\mathfrak K}_{\mathbf A}$ follows basically from \cite[Lemma 2.5]{CibilsMarcos}.
Here is the adaption of the proof.
By Lemma~\ref{homsforkleisli} we have
$${{\mathfrak K}_{\mathbf A}}(X,Y)={\mathcal H}(X,AY)=
{\mathcal H}(X,\bigoplus_{i\in I}E_iY).$$
Hence, for $X=\widehat E_iM$ and $Y=M$ we get that ${{\mathfrak K}_{\mathbf A}}(X,Y)$ contains
the identity on $\widehat E_i M$, which is clearly invertible. Therefore $M\simeq \widehat E_iM$
in ${\mathfrak K}_{\mathbf A}$.
\dickebox

\begin{Cor}\label{otheradjointpairs}
If $(S,T)$ is another adjoint pair inducing this
monad ${\mathbf A}$ (where $S:{\mathcal H}\lra{\mathcal G}$),
then let $K:{\mathfrak K}_{\mathbf A}\lra{\mathcal G}$ be the functor
from Section~\ref{KleisliEMsection}.
If now there are self-equivalences $E_i^0$ of $\mathcal G$
such that $K\widehat E_i=E_i^0K$ for all $K$, then also $E_i^0V\simeq V$ for all $V$
in the image of $S$.
\end{Cor}

Proof. Recall that by Section~\ref{KleisliEMsection}
the Kleisli category is initial amongst all adjoint pairs inducing
a fixed monad $\mathbf A$ in the sense that the functor $K$ exists and satisfies $KS_K=S$.
Then by Proposition~\ref{widehatEiMisisotoM} we have that
$\widehat E_iV\simeq V$ for all $V$ in ${\mathfrak K}_{\mathbf A}$. Recall that all objects in
${\mathfrak K}_{\mathbf A}$ are in the image of $S_K$. Hence $V=S_KW$ and
$$E_i^0SW=E_i^0KS_KW=K\widehat E_iS_KW=K\widehat E_iV\simeq KV=KS_KW=SW.$$
This shows the Corollary. \dickebox.

\begin{Rem}
Note that if the counit $\epsilon:S_KT_K\lra \textup{id}_{{\mathfrak K}_{\mathbf A}}$
(or actually even more generally
the counit $\epsilon:ST\lra \textup{id}_{\mathcal G}$) splits functorially, then
by Proposition~\ref{monadsadmittingsection} we have an equivalence
${\mathcal G}\simeq {\mathbf A}-\textup{Mod}_{\mathcal H}$ up to direct factors,
an equivalence in case the categories are idempotent complete, and therefore we may
assume that this equivalence is actually an equality.
By the universal property of the Kleisli category, we also get that
$${\mathcal G}\simeq {\mathbf A}-\textup{Mod}_{\mathcal H}\simeq {\mathfrak K}_{\mathbf A}$$
if the categories are idempotent complete.
\end{Rem}

\begin{Theorem}\label{Resultincaseoftrivialinertiagroup}
Let $\mathcal H$ and $\mathcal G$ be additive Krull-Schmidt categories. Let
$$
\xymatrix{{\mathcal H}\ar@/^/[r]^S&{\mathcal G}\ar@/^/[l]^T}
$$
be functors giving a normal subgroup situation $TS=\bigoplus_{i\in I}E_i$,
such that $(S,T)$ is an adjoint pair, and such that $E_{i_0}=\textup{id}_{\mathcal H}$.
Suppose that for any $i\in I$ there is a self-equivalence $\widehat E_i$ of $\mathcal G$ with
$SE_i=\widehat E_iS$ and $T\widehat E_i=E_iT$.
Then for any indecomposable object $W$ of $\mathcal H$ with $E_jW\simeq W$ if and only if $j=i_0$, we get
$S(W)$ is indecomposable as well.
\end{Theorem}

Proof. If $S(W)=V_1\oplus V_2$, then $TSW=TV_1\oplus TV_2$ and since
$\textup{id}_{\mathcal H}$ is a direct factor of $TS$, $W$
is a direct factor of
either $TV_1$ or $TV_2$, w.l.o.g. of $TV_1$, and we may assume that $V_1$
is indecomposable. Since for any indecomposable object $W$ of $\mathcal H$
with $E_jW\simeq W$ if and only if $j=i_0$, we have
$E_iW\not\simeq W$ for all $i\in I\setminus\{i_0\}$.
Since $E_j$ is a self-equivalence, $E_jW$ is indecomposable again
for all $j\in I$.
But then $E_jW$ is a direct factor of $E_jTV_1=T\widehat E_jV_1$.
However, $\widehat E_jV_1\simeq V_1$ by Proposition~\ref{ehatmapstoisomorphiccopies}
and Corollary~\ref{otheradjointpairs} and hence
$E_jW$ is a direct factor of $TV_1$ for all $j\in I$.
Therefore $AW=TSW$ is a direct factor of $TV_1$. But $TSW=TV_1\oplus TV_2$.
Hence $TV_2=0$ by the Krull-Schmidt property of $\mathcal H$.
But this implies
$$0={\mathcal H}(W,TV_2)\simeq {\mathcal G}(SW,V_2)={\mathcal G}(V_1\oplus V_2,V_2)$$
which has a direct factor ${\mathcal G}(V_2,V_2)$, containing the identity on $V_2$.
Therefore $V_2=0$ and $S(W)$ is indecomposable. \dickebox

\begin{Rem}
Note that for $\mathcal G$ being the Kleisli category or the Eilenberg-Moore category of the adjunction, then Proposition~\ref{EiliftstoG} shows the existence of $\widehat E_i$ for all $i$.
\end{Rem}

\section{The inertia category}

\label{generalcasesection}

We can now imitate the construction of an inertia group.

\subsection{On orbit categories and Kleisli categories}

Let $(A,\mu,\eta)$ be a monad on $\mathcal H$, and
suppose $A=\bigoplus_{i\in I} E_i$ for some automorphisms $E_i$ of $\mathcal H$.
We suppose that $\{E_i\;|\;i\in I\}$ with the multiplication $\mu$
forms a group $\Gamma$ of automorphisms of $\mathcal H$.
If $\mu:A\circ A\lra A$ is the structure map, then
writing $A=\bigoplus_{i\in I}E_i$ we get by restriction natural
transformations $\mu^{i,j}:E_i\circ E_j\lra E_k$, and $\mu^{i,j}=0$ for almost all $k$,
and for a unique $k=k(i,j)$ this is an equivalence.

Let now $(S,T)$ be a pair of adjoint functors
$$
\xymatrix{{\mathcal H}\ar@/^/[r]^S&{\mathcal G}\ar@/^/[l]^T}
$$
inducing the monad $A$ giving a normal subgroup situation $A=\bigoplus_{i\in I}E_i$
such that the structure map on $A$ induces a group law on the $E_i$, as indicated above.
Then the hypothesis $E_iA\simeq AE_i$ for all $i$ from Proposition~\ref{EiliftstoG} holds true.

We recall the construction of the orbit category as is displayed in \cite{CibilsMarcos}.
We use the notation introduced there, though we mention that in Keller's
\cite{triangulatedorbitcategories} the notion ${\mathcal C}/\Gamma$ is used for the category
which is denoted by ${\mathcal C}[\Gamma]$ in \cite{CibilsMarcos}.

Let $\mathcal C$ be a category admitting arbitrary direct sums and let $\Gamma$ be
a group of self-equivalences of $\mathcal C$. Then we may form the orbit category
${\mathcal C}[\Gamma]$.
The objects of $\mathcal C$ and of ${\mathcal C}[\Gamma]$ coincide.
Further,
$${{\mathcal C}[\Gamma]}(X,Y)=\bigoplus_{g\in\Gamma}{\mathcal C}(X,gY).$$
Note that we have
$$
{{\mathcal C}[\Gamma]}(X,Y)=
\left(\prod_{(g_1,g_2)\in G\times G}{\mathcal C}(g_1X,g_2Y)\right)^G
$$
where $G$ acts diagonally. This observation can be found in e.g.
Asashiba~\cite{AsashibaGalois1}. Composition of morphisms is then clear with the above formula.

\medskip

The orbit category of a group is actually a special case of a Kleisli category.

\begin{Cor}\label{Kleisliandorbitcat} \label{RGammaisrightadjoint}
Let ${\mathbf A}=(A,\mu,\eta)$ be a monad,
suppose $A=\bigoplus_{i\in I} E_i$ for self-equivalences $E_j$ of $\mathcal H$, and
suppose that $\{E_i\;|\;i\in I\}$ with the structure map $\mu$ forms a group $\Gamma$
of self-equivalences of $\mathcal H$. Then,
the orbit category ${\mathcal H}[\Gamma]$ is naturally isomorphic to
${\mathfrak K}_{\mathbf A}$, the Kleisli category and $\mathbf A$ is the monad corresponding to
an adjoint pair $(N[\Gamma],R[\Gamma])$
$$
\xymatrix{{\mathcal C}\ar@/^/[r]^{N[\Gamma]}&{\mathcal C}[\Gamma]\ar@/^/[l]^{R[\Gamma]}}
$$
\end{Cor}

Proof.
This is a direct consequence of Lemma~\ref{homsforkleisli} where we showed that
$${\mathcal H}[\Gamma]\simeq {\mathfrak K}_{\mathbf A}.$$
Moreover, $N[\Gamma]$ identifies with $S_{{\mathfrak K}_{\mathbf A}}$ and
$R[\Gamma]$ identifies with $T_{{\mathfrak K}_{\mathbf A}}$.
\dickebox

\medskip

We recall a concept from Asashiba~\cite{AsashibaGalois1}.

\begin{Def}\cite{AsashibaGalois1}\label{gammainvariantfunctordef}
Let $\mathcal C$ be a category with an action of a group $\Gamma$ in the sense that
there is a set of self-equivalences $E_g$ of $\mathcal C$ which form a group $\Gamma$.
We denote by $E_1$ the identity self-equivalence.

Then for a functor $F:{\mathcal C}\lra{\mathcal D}$ a collection of natural transformations
$\phi_g:F\lra FE_g$ for $E_g\in \Gamma$ is called $\Gamma$-adjuster if
$$\xymatrix{F\ar[r]^-{\phi_g}\ar[dr]_-{\phi_{hg}}&FE_g\ar[d]^-{\phi_hE_g}\\ &FE_{hg}\ar@{=}[r]&FE_hE_g}$$
is commutative.
A functor $F$ for which there is a $\Gamma$-adjuster is called a $\Gamma$-invariant functor.
\end{Def}

Asashiba remarks that then $\phi_1=\textup{id}_F$ automatically.

\begin{Lemma}\cite[Proposition 2.6]{AsashibaGalois1}\label{Asashibainvariantfunctorlemma}
Let $\mathcal C$ be a category with an action of a group $\Gamma$ in the sense that
there is a set of self-equivalences $E_g$ of $\mathcal C$ which form a group $\Gamma$. Then
\begin{enumerate}
\item \label{Asashibainvariantfunctorlemmai}
The natural functor $N[\Gamma]:{\mathcal C}\lra {\mathcal C}[\Gamma]$ is $\Gamma$-invariant with adjuster $(\nu_g)_{g\in\Gamma}$, where for each $g\in\Gamma$ and object $X$ of $\mathcal C$ we define $$\nu_{g,X}:=(\delta_{h_1,h_2g}\textup{id}_{h_1X})_{(h_1,h_2)\in\Gamma\times\Gamma}\in Hom_{{\mathcal C}[\Gamma]}(N[\Gamma]X,N[\Gamma]gX).$$
\item \label{Asashibainvariantfunctorlemmaii}
For each $\Gamma$-invariant functor ${\mathcal C}\stackrel{(F,\phi)}\lra{\mathcal D}$
there is a unique functor ${\mathcal C}[\Gamma]\stackrel{H}\lra{\mathcal D}$
such that $(F,\phi)= (HN[\Gamma],H\nu)$, i.e. in particular the diagram
$$\xymatrix{{\mathcal C}\ar[r]^{F}\ar[d]_{N[\Gamma]}&{\mathcal D}\\
{\mathcal C}[\Gamma]\ar[ru]_H
}$$
is commutative.
\item \label{Asashibainvariantfunctorlemmaiii}
For each $\Gamma$-invariant functor ${\mathcal C}\stackrel{(F,\phi)}\lra{\mathcal D}$
there is an up to isomorphism unique functor ${\mathcal C}[\Gamma]\stackrel{H}\lra{\mathcal D}$
such that $(F,\phi)\simeq (HN[\Gamma],H\nu)$,
\end{enumerate}
\end{Lemma}

We have to deal with submonads and orbit categories with respect to a subgroup.

\begin{Prop}\label{factoringorbitcategories}
Let $\mathcal H$ be an additive category, let $\Gamma$ be a group of automorphisms
and let $\Gamma_0$ be a subgroup of $\Gamma$. Consider the associated adjoint
pairs
$$
\xymatrix{{\mathcal H}\ar@/^/[r]^{S_\Gamma}&{\mathcal H}[\Gamma]\ar@/^/[l]^{T_\Gamma}}
$$
and
$$
\xymatrix{{\mathcal H}\ar@/^/[r]^{S_{\Gamma_0}}&{\mathcal H}[\Gamma_0]\ar@/^/[l]^{T_{\Gamma_0}}}
$$
given by the Kleisli categories associated to the monads
$A=\bigoplus_{\gamma\in\Gamma}\gamma$
respectively $A_0=\bigoplus_{\gamma\in\Gamma_0}\gamma$.
Then there is an adjoint pair
$(S_{\Gamma^0},T_{\Gamma^0})$
$$
\xymatrix{{\mathcal H}[\Gamma_0]\ar@/^/[r]^{S_{\Gamma^0}}&{\mathcal H}[\Gamma]\ar@/^/[l]^{T_{\Gamma^0}}}
$$
such that
$$S_\Gamma=S_{\Gamma^0}\circ S_{\Gamma_0}\textup{ and }
T_\Gamma=T_{\Gamma_0}\circ T_{\Gamma^0}.$$
\end{Prop}

Proof.
By Lemma~\ref{Asashibainvariantfunctorlemma}.\ref{Asashibainvariantfunctorlemmai}
the functor
$S_\Gamma:{\mathcal H}\lra {\mathcal H}[\Gamma]$ is $\Gamma$-invariant in
the sense of Definition~\ref{gammainvariantfunctordef}. Hence, by
restriction to the subgroup $\Gamma_0$, $S_\Gamma$ is $\Gamma_0$-invariant as well. By
Lemma~\ref{Asashibainvariantfunctorlemma}.\ref{Asashibainvariantfunctorlemmaii}
factorizes through ${\mathcal H}[\Gamma_0]$. This shows the existence of $S_{\Gamma^0}$.

We give the functor explicitly.
$${\mathcal H}[\Gamma_0]\stackrel{S_{\Gamma^0}}\lra {\mathcal H}[\Gamma]$$
satisfying
$S_\Gamma=S_{\Gamma^0}\circ S_{\Gamma_0}$.
We put $S_{\Gamma^0}(X)=X$ on objects. Then, since
$S_{\Gamma_0}(X)=X$ on objects and $S_{\Gamma}(X)=X$ on objects
we have $S_\Gamma=S_{\Gamma^0}\circ S_{\Gamma_0}$ on objects.
On morphisms we have $S_\Gamma(f)=\bigoplus_{\gamma\in\Gamma}\gamma f$
and $S_{\Gamma_0}(f)=\bigoplus_{\gamma\in\Gamma_0}\gamma f$
for any $f\in Hom_{\mathcal H}(X,Y)$.
Then
$$
{{\mathcal H}[\Gamma_0]}(X,Y)
={\mathcal H}(X,\bigoplus_{\gamma\in\Gamma_0}\gamma Y)
$$
and for $f\in {\mathcal H}(X,\bigoplus_{\gamma\in\Gamma_0}\gamma Y)$
we put
$$S_{\Gamma^0}(f):=\bigoplus_{\Gamma_0\sigma\in\Gamma/\Gamma_0}\sigma f.$$
Then by definition $S_\Gamma=S_{\Gamma^0}\circ S_{\Gamma_0}$ on morphisms.

Similarly, $T_{\Gamma_0}(X)=\bigoplus_{\Gamma_0\sigma\in\Gamma/\Gamma_0}\sigma X$
and $T_{\Gamma_0}(f)=f$ for each morphism $f$.
We need to explain briefly the meaning here.
Note that
$${\mathcal H}[\Gamma_0](X,Y)={\mathcal H}(X,\bigoplus_{\gamma\in\Gamma_0}\gamma Y)
={\mathcal H}[\Gamma_0](X,Y)={\mathcal H}(\bigoplus_{\gamma\in\Gamma_0}\gamma X,\bigoplus_{\gamma\in\Gamma_0}\gamma Y)^{\Gamma_0}$$ and let
$$f\in{\mathcal H}[\Gamma](X,Y)={\mathcal H}(X,\bigoplus_{\gamma\in\Gamma}\gamma Y)={\mathcal H}(\bigoplus_{\gamma\in\Gamma}\gamma X,\bigoplus_{\gamma\in\Gamma}\gamma Y)^{\Gamma}.$$ Then
\begin{eqnarray*}
T_{\Gamma_0}(f)\in {\mathcal H}[\Gamma_0](T_{\Gamma_0}X,T_{\Gamma_0}Y)&=&
{\mathcal H}[\Gamma_0](\bigoplus_{\Gamma_0\sigma\in\Gamma/\Gamma_0}\sigma X,\bigoplus_{\Gamma_0\sigma\in\Gamma/\Gamma_0}\sigma Y)\\
&=&
{\mathcal H}(\bigoplus_{\gamma\in\Gamma_0}\gamma\bigoplus_{\Gamma_0\sigma\in\Gamma/\Gamma_0}\sigma X,\bigoplus_{\gamma\in\Gamma_0}\gamma\bigoplus_{\Gamma_0\sigma\in\Gamma/\Gamma_0}\sigma Y)^{\Gamma_0}\\
&=&
{\mathcal H}(\bigoplus_{\gamma\in\Gamma}\gamma X,\bigoplus_{\gamma\in\Gamma}\gamma Y)^{\Gamma_0}
\end{eqnarray*}
and $T_{\Gamma_0}(f)$ just maps $f$ in the space of $\Gamma$-fixed points to the space of $\Gamma_0$-fixed points. This is clearly a functor.
Furthermore, since $(S,T)$ is an adjoint pair, also
$(S_{\Gamma^0},T_{\Gamma^0})$ is an adjoint pair.

This shows the statement. \dickebox

\begin{Rem}
Note that the unit of the adjunction $(S_{\Gamma^0},T_{\Gamma_0})$ leads to some
Mackey formula involving double classes.
\end{Rem}

\begin{Rem}
In~\cite{Chavez}, in case $\Gamma$ is cyclic, Ch\'avez gave
a criterion when the orbit category is idempotent complete (i.e. all idempotents split). In general
the orbit category is not idempotent complete.
Indeed,
let $A=Mat_2(K)$ for some field $K$. Conjugation by $M:=\left(\begin{array}{cc}0&1\\1&0\end{array}\right)$
is an automorphism of order $2$ of $A$, yielding an action of the cyclic group of order $2$ on $A-mod$.
Then, for a simple $A$-module $S$ we get
$$Hom_{A-mod[C_2]}(S,S)=Hom_A(S,S\oplus\ ^MS)=K\times K$$
has non trivial idempotents, whereas $S$ is indecomposable. Moreover, in
Krull-Schmidt categories endomorphism algebras of indecomposable objects are local.
The above endomorphism algebra is not local.

Similarly, let $A$ be the Kronecker algebra and $\sigma$ be the automorphism swapping the two arrows
of the Kronecker quiver. This induces an action of the cyclic group of order $2$ on the module category.
The simples are fixed under the automorphism.
\end{Rem}

\begin{Rem}\label{Karoubienvelope}
Recall the Karoubi envelope of a category.
Let $\mathcal C$ be a category. Then $Kar({\mathcal C})$ has a objects
the pairs $(C,e)$ for objects $C$ of $\mathcal C$ and
$e^2=e\in End_{\mathcal C}(C)$.
Further,
$$Kar({\mathcal C})((C,e),(D,f))=\{\alpha\in {\mathcal C}(C,D)\;|\;
\alpha\circ e=f\circ\alpha\}$$
Note that there is a fully faithful embedding
$${\mathcal C}\hookrightarrow Kar({\mathcal C})$$
given by mapping an object $C$ to $(C,1_C)$ and a morphism to this same morphism.
It is well-known, and not difficult to show, that if $\mathcal C$ is idempotent complete,
then the natural embedding ${\mathcal C}\hookrightarrow Kar({\mathcal C})$
is an equivalence.

If $S:{\mathcal H}\lra{\mathcal G}$ is a functor. Then we define a functor
$$\widehat S:Kar({\mathcal H})\lra Kar({\mathcal G})$$
by $\widehat S(X,e):=(S(X),S(e))$ and $\widehat S(\beta):=S(\beta)$ for any morphism $\beta$.
Since $$S(\beta)S(e)=S(\beta e)=S(f\beta)=S(f)S(\beta)$$
with the obvious notations, this is indeed well-defined.
\end{Rem}

\begin{Rem}
We should mention that an additive category, whose morphism spaces are
$R$-modules of finite length over some commutative ring $R$, is Krull-Schmidt
if and only if it has split idempotents (i.e. is idempotent complete).
This was proved most recently  by Amit Shah \cite{AmitShah}, and may have
been known earlier to experts.
\end{Rem}

\begin{Lemma}\label{adjointpairsonKaroubienvelope}
Let  now $(S,T)$ be a pair of adjoint functors
$$
\xymatrix{{\mathcal H}\ar@/^/[r]^S&{\mathcal G}\ar@/^/[l]^T}
$$
with unit $\eta:\textup{id}_{\mathcal H}\lra TS$ and counit $\epsilon:ST\lra \textup{id}_{\mathcal G}$.
Then
$$
\xymatrix{Kar({\mathcal H})\ar@/^/[r]^{\widehat S}&Kar({\mathcal G})\ar@/^/[l]^{\widehat T}}
$$
is again an adjoint pair with counit $\widehat\epsilon_{(X,e)}:=\epsilon_X$ and
unit $\widehat\eta_{(X,e)}:=\eta_X$.
\end{Lemma}

Proof.
Indeed, we first observe by direct inspection
that $\widehat\epsilon$ and $\widehat\eta$ are natural
transformations.
Since $(S,T)$ is an adjoint pair, the compositions
$$S\stackrel{S\eta}{\lra}STS\stackrel{\epsilon S}{\lra}S$$
and
$$T\stackrel{\eta T}{\lra}TST\stackrel{T\epsilon}{\lra}T$$
are the identity.
But then,
$$\widehat S\stackrel{\widehat S\widehat\eta}{\lra}\widehat S\widehat T\widehat S\stackrel{\widehat\epsilon \widehat S}{\lra}\widehat S$$
and
$$\widehat T\stackrel{\widehat \eta \widehat T}{\lra}\widehat T\widehat S\widehat T\stackrel{\widehat T\widehat\epsilon}{\lra}\widehat T$$
are the identity as well since the natural transformations on the Karoubi envelope are defined just as on the original categories. Hence, by \cite[Chapter IV, Section 1, Theorem 2]{Maclane}
$(\widehat S,\widehat T)$ is an adjoint pair with unit $\widehat\eta$ and counit $\widehat\epsilon$. \dickebox

\medskip

\subsection{Preparing adjoint functors for the Inertia category}

\begin{Lemma}\label{monadonorbitcategory}
Let $(A,\mu,\eta)$ be a monad giving a normal subgroup situation,
i.e. $A=\bigoplus_{i\in I}E_i$ for automorphisms $E_i$, suppose
that $\{E_i\;|\;i\in I\}$ with the structure map $\mu$
forms a group $\Gamma$ of automorphisms.
Let $\Gamma_0$ be a subgroup of $\Gamma$  of index $|\Gamma:\Gamma_0|$.
Suppose that
$$
\xymatrix{{\mathcal H}\ar@/^/[r]^S&{\mathcal H}[\Gamma]\ar@/^/[l]^T}
$$
is an adjoint pair of additive functors realising $A$,
and suppose that the counit of the adjunction $(S,T)$ is locally split.
Then by Proposition~\ref{factoringorbitcategories} the adjoint pair $(S,T)$
induces naturally an adjoint pair $(S[\Gamma^0],T[\Gamma^0])$
$$
\xymatrix{{\mathcal H}[\Gamma_0]\ar@/^/[r]^-{S[\Gamma^0]}&
{\mathcal H}[\Gamma]\ar@/^/[l]^-{T[\Gamma^0]}}
$$
satisfying $S=S[\Gamma^0]\circ S[\Gamma_0]$ and $T=T[\Gamma_0]\circ T[\Gamma^0]$.
Moreover, each $E_i$ induces an automorphism $E_i^0$ of ${\mathcal H}[\Gamma_0]$
with $E_i^0N[\Gamma_0]=N[\Gamma_0]E_i$ and $E_iR[\Gamma_0]=R[\Gamma_0]E_i^0$.
Furthermore, on objects we get
$$T[\Gamma^0]\circ S[\Gamma^0]=\bigoplus_{E_i\Gamma_0\in\Gamma/\Gamma_0}E_i^0.$$
\end{Lemma}

Proof.
By construction
$${{\mathcal H}[\Gamma_0]}(X,Y)={\mathcal H}(X,\bigoplus_{g\in\Gamma_0}gY)$$
and since the counit is locally split, by Proposition~\ref{EiliftstoG}
each $E_i$ induces an automorphism $E_i^0$ of ${\mathcal H}[\Gamma_0]$
with $E_i^0N[\Gamma_0]=N[\Gamma_0]E_i$ and $E_iR[\Gamma_0]=R[\Gamma_0]E_i^0$.
Similarly,
each $E_i$ induces an automorphism $\widehat E_i$ of
${\mathcal H}[\Gamma]$
with $\widehat E_iS=SE_i$ and $E_iT=T \widehat E_i$.

\begin{eqnarray*}
{{\mathcal H}[\Gamma_0]}(-,T[\Gamma^0]S[\Gamma^0]Y)
&=&{{\mathcal H}[\Gamma]}(S[\Gamma^0]-,S[\Gamma^0]Y)\\
&=&{{\mathcal H}}(-,\bigoplus_{E_i\in\Gamma}E_iY)\\
&=&{{\mathcal H}}(-,\bigoplus_{E_i\Gamma_0\in\Gamma/\Gamma_0}E_i\bigoplus_{E_j\in\Gamma_0}E_jY)\\
&=&{{\mathcal H}}(\bigoplus_{E_i\Gamma_0\in\Gamma/\Gamma_0}E_i^{-1}-,\bigoplus_{E_j\in\Gamma_0}E_jY)\\
&=&{{\mathcal H}[\Gamma_0]}(\bigoplus_{E_i\Gamma_0\in\Gamma/\Gamma_0}{(E_i^0)}^{-1}-,Y)\\
&=&{{\mathcal H}[\Gamma_0]}(-,\bigoplus_{E_i\Gamma_0\in\Gamma/\Gamma_0}E_i^0Y)
\end{eqnarray*}

By Yoneda's lemma we have
$$\left(T[\Gamma^0]S[\Gamma^0]\right)\simeq \bigoplus_{E_i\Gamma_0\in\Gamma/\Gamma_0}E_i^0$$
on objects.
\dickebox

\subsection{Clifford's theorem for orbit categories, the main result}

We are now ready to prove our main result.

\begin{Theorem}\label{CliffordforMonads}
Let $R$ be a commutative Noetherian ring.
Suppose that
$$
\xymatrix{{\mathcal H}\ar@/^/[r]^S&{\mathcal G}\ar@/^/[l]^T}
$$
is an adjoint pair of $R$-linear functors between $R$-linear idempotent
complete categorie. Assume that that $Hom$-spaces
are of finite $R$-length in ${\mathcal H}$ and in ${\mathcal G}$, and assume that $(S,T)$
give a normal subgroup situation.
Let $TS=\bigoplus_{i\in I}E_i$ for automorphisms $E_i$, suppose
that $\Gamma:=\{E_i\;|\;i\in I\}$ forms a group of automorphisms
and for an indecomposable object $M$
let $I_S(M):=\{i\in I\;|\;E_iM\simeq M\},$ and let
$\Gamma_M:=\{E_i\;|\;i\in I_{S}(M)\}$.
Suppose that $\mathcal G$ is  the Karoubi envelope of the
Kleisli construction of the monad
(or, what is the same, the orbit category with respect to $\Gamma$),
and suppose that $\Gamma$ is finite.
Then we have a commutative diagram of adjoint pairs
$$\xymatrixcolsep{5pc}\xymatrixrowsep{5pc}
\xymatrix{
{\mathcal H}\ar@/^/[rr]^S\ar@/^1pc/[dr]|-{S[\Gamma_M]}&&{\mathcal G}\ar@/^/[ll]^T\ar@/^1pc/[dl]|-{T[\Gamma_M^0]}\\
&Kar({\mathcal H}[\Gamma_M])\ar@/^1pc/[ur]|-{S[\Gamma_M^0]}\ar@/^1pc/[ul]|-{T[\Gamma_M]}
}
$$
and for any indecomposable direct
factor $M_0$ of $S[\Gamma_M](M)$ we have that
$S[\Gamma_M^0](M_0)$ is indecomposable.
\end{Theorem}

Proof. We first suppose that $\mathcal G$ is actually the Kleisli
construction.
Since $\Gamma$ is a group, $E_iA\simeq AE_i$ for all $i$ and the isomorphisms
satisfy the $1$-cocycle condition.
If $\mathcal G$ is the Kleisli construction of the monad, then, by
Proposition~\ref{EiliftstoG} and Remark~\ref{liftingthegroup},
we get that $\Gamma_M$ acts on
$\mathcal G$ as group of automorphisms $\widehat\Gamma_M$
with automorphisms $\widehat E_i^M$ of $\mathcal G$.

Moreover, by same argument, there are automorphisms
$E_i^M$ of ${\mathcal H}[\Gamma_M]$ forming a group $\Gamma_M$ with
$$E_i^MS[\Gamma_M]=S[\Gamma_M]E_i\textup{ and }
T[\Gamma_M]E_i^M=E_iT[\Gamma_M].$$

By Lemma~\ref{RGammaisrightadjoint} we get that
$(S[\Gamma_M],T[\Gamma_M])$ is
an adjoint pair
$$
\xymatrix{
{\mathcal H}\ar@/^/[r]^-{S[\Gamma_M]}&
{\mathcal H}[\Gamma_M]\ar@/^/[l]^{T[\widehat\Gamma_M]}
}
$$

Since $\Gamma$ is finite, $\Gamma_M$ is finite as well. Further, since all
homomorphisms are $R$-modules of finite length, the same holds true for
the orbit category with respect to $\Gamma_M$ and to $\Gamma$. Since $R$ is Noetherian,
this is still true for the Karoubi envelopes. Hence, by \cite{AmitShah}
the Karoubi envelopes
of the orbit categories with respect to $\Gamma_M$ and to $\Gamma$ are Krull-Schmidt categories.
By Lemma~\ref{monadonorbitcategory} we have a commutative diagram of adjoints
$$\xymatrixcolsep{5pc}\xymatrixrowsep{5pc}
\xymatrix{
{\mathcal H}\ar@/^/[rr]^S\ar@/^1pc/[dr]|-{S[\Gamma_M]}&&
{\mathcal G}\ar@/^/[ll]^T\ar@/^1pc/[dl]|-{T[\Gamma_M^0]}\\
&{\mathcal H}[\Gamma_M]\ar@/^1pc/[ur]|-{S[\Gamma_M^0]}\ar@/^1pc/[ul]|-{T[\Gamma_M]}
}
$$

Now, in case $\mathcal G$ is the Karoubi envelope of the Kleisli construction,
since $\mathcal H$ and $\mathcal G$ are idempotent complete, we can
replace ${\mathcal H}[\Gamma_M]$ by its Karoubi envelope
$\widehat{{\mathcal H}[\Gamma_M]}$.
Then, using Lemma~\ref{adjointpairsonKaroubienvelope}, in the diagram
$$\xymatrixcolsep{5pc}\xymatrixrowsep{5pc}
\xymatrix{
{\mathcal H}\ar@/^/[rr]^S\ar@/^1pc/[dr]|-{S[\Gamma_M]}&&
{\mathcal H}[\Gamma]\ar@/^/[ll]^T\ar@/^1pc/[dl]|-{T[\Gamma_M^0]}\ar@{^{(}->}[r]|-{D}&{\mathcal G}\\
&{\mathcal H}[\Gamma_M]\ar@/^1pc/[ur]|-{S[\Gamma_M^0]}
\ar@/^1pc/[ul]|-{T[\Gamma_M]}\ar@{^{(}->}[d]|-{C}\\
&\widehat{{\mathcal H}[\Gamma_M]}
}
$$
the functor $T[\Gamma_M]$ and $S[\Gamma_M^0]$ lift to
$\widehat{{\mathcal H}[\Gamma_M]}$ such that the diagram
$$\xymatrixcolsep{5pc}\xymatrixrowsep{5pc}
\xymatrix{
{\mathcal H}\ar@/^/[rr]^S\ar@/^1pc/[dr]|-{S[\Gamma_M]}&&
{\mathcal H}[\Gamma]\ar@/^/[ll]^T\ar@/^1pc/[dl]|-{T[\Gamma_M^0]}\ar@{^{(}->}[r]|-{D}&{\mathcal G}\ar@/_1pc/[lldd]|-{\widehat{T[\Gamma_M^0]}}\\
&{\mathcal H}[\Gamma_M]\ar@/^1pc/[ur]|-{S[\Gamma_M^0]}
\ar@/^1pc/[ul]|-{T[\Gamma_M]}\ar@{^{(}->}[d]|-{C}\\
&\widehat{{\mathcal H}[\Gamma_M]}\ar@/^1pc/[luu]|-{\widehat{T[\Gamma_M]}}
\ar@/_1pc/[rruu]|-{\widehat{S[\Gamma_M^0]}}
}
$$
is commutative in the natural sense. Moreover, $(\widehat{S[\Gamma_M]},\widehat{T[\Gamma_M]})$ and
$(\widehat{S[\Gamma_M^0]},\widehat{T[\Gamma_M^0]})$ are adjoint pairs again.
Furthermore, $E_i^M$ extend to automorphisms of $\widehat{{\mathcal H}[\Gamma_M]}$,
and each $\widehat E_i^M$ extend to an automorphism of $\mathcal G$,
and in order to avoid an additional notational burden, we denote the extension of
$E_i^M$ to the Karoubi envelope by $E_i^M$ again, and likewise for $\widehat E_i^M$.

Let
$$CS[\Gamma_M](M)=\bigoplus_{j\in J}M_j$$
for indecomposable objects $M_j$ of $\widehat{{\mathcal H}[\Gamma_M]}$.
Then, using Lemma~\ref{monadonorbitcategory},
$$
\widehat{T[\Gamma_M^0]}\widehat{S[\Gamma_M^0]}(M_j)=\left(
\bigoplus_{E_i\in\Gamma/\Gamma_M}E_i^M(M_j)
\right).
$$

Let $M_j$ be an indecomposable direct factor of $CS[\Gamma_M](M)$.

By definition for each representative $E_i$ of a class different from
$\Gamma_M$ in $\Gamma/\Gamma_M$
we have $E_iM\not\simeq M$. By the definition of $\Gamma_M$,
$$\widehat{T[\Gamma_M]}CS[\Gamma_M]M\simeq
T[\Gamma_M]S[\Gamma_M]M\simeq  \bigoplus_{E_i\in\Gamma_M}E_iM\simeq \bigoplus_{E_i\in\Gamma_M}M.$$

The object $M$ being indecomposable, $M_j$ being a direct factor of
$CS[\Gamma_M]M$, hence $\widehat{T[\Gamma_M]}M_j$ is a direct factor of
$\bigoplus_{E_i\in\Gamma_M}M.$
Therefore, by the Krull-Schmidt theorem on $\mathcal H$,
$\widehat{T[\Gamma_M]}M_j$ is a direct
sum of,  $n_j$ say, copies of $M$. In particular $E_i^MM_j\not\simeq M_j$
whenever $E_i^M\not\in\Gamma_M$, and moreover, since
$$M^{|\Gamma_M|}=T[\Gamma_M]S[\Gamma_M]M=\widehat{T[\Gamma_M]}CS[\Gamma_M]M=
\widehat{T[\Gamma_M]}(\bigoplus_{j\in J}M_j)=
\bigoplus_{j\in J}\widehat{T[\Gamma_M]}M_j=\bigoplus_{j\in J}M^{n_j}$$
we get $\sum_{j\in J}n_j=|\Gamma_M|$.

Fix now $j\in J$ and consider the indecomposable
direct factor $M_j$ of $CS[\Gamma_M](M)$.
Since $M_j$ is a direct factor of $(\widehat{T[\Gamma_M^0]}\widehat{S[\Gamma_M^0]})(M_j)$,
there is an indecomposable direct factor $X_j$ of
$\widehat{S[\Gamma_M^0]}(M_j)=X_j\oplus Y_j$ in $\mathcal G$ such that
$M_j$ is a direct factor of $\widehat{T[\Gamma_M^0]}(X_j)$ and
$M^{n_j}$ is a direct factor of $\widehat{T[\Gamma_M]}\widehat{T[\Gamma_M^0]}(X_j)$.

Recall that  there are automorphisms $\widehat E_i^M$ of $\mathcal G$
such that
$$
E_i^M\widehat{T[\Gamma_M^0]}=\widehat{T[\Gamma_M^0]}\widehat E_i^M
$$
and by
Proposition~\ref{ehatmapstoisomorphiccopies} we have
$\widehat E_i^M(V)\simeq V$ for all $E_i\in\Gamma$ and objects $V$.
Since $M_j$ is a direct factor of $\widehat{T[\Gamma_M^0]}X_j$, we get that
$E_i^M(M_j)$ is a direct factor of
$$E_i^M\widehat{T[\Gamma_M^0]}X_j=\widehat{T[\Gamma_M^0]}\widehat E_i^MX_j\simeq \widehat{T[\Gamma_M^0]}X_j$$
for all $E_i\in\Gamma$.
Since for each non trivial class $E_i^M\Gamma_M$ of $\Gamma/\Gamma_M$ we have
$E_i^M(M_j)\not\simeq M_j$, as seen above,
$\bigoplus_{E_i\Gamma_M\in\Gamma/\Gamma_M}E_i^M(M_j)$ is a
direct factor of $\widehat{T[\Gamma_M^0]}(X_j)$, using the Krull-Schmidt property
on $\widehat{{\mathcal H}[\Gamma_M]}$. Furthermore
$$\widehat{T[\Gamma_M]}\left(\bigoplus_{E_i\Gamma_M\in\Gamma/\Gamma_M}E_i^M(M_j)\right)=
\bigoplus_{E_i\Gamma_M\in\Gamma/\Gamma_M}E_i\widehat{T[\Gamma_M]}(M_j)=
\bigoplus_{E_i\Gamma_M\in\Gamma/\Gamma_M}E_iM^{n_j}
$$
is a direct factor of $\widehat{T[\Gamma_M]}\widehat{T[\Gamma_M^0]}(X_j)$.

However,
\begin{eqnarray*}
\widehat{T[\Gamma_M]}\widehat{T[\Gamma_M^0]}\widehat{S[\Gamma_M^0]}M_j
&=&\widehat{T[\Gamma_M]}
\left(\bigoplus_{E_i^M\Gamma_M\in\Gamma/\Gamma_M}E_i^M\right)M_j\\
&=&\left(\bigoplus_{E_i^M\Gamma_M\in\Gamma/\Gamma_M}E_i\right)\widehat{T[\Gamma_M]}M_j\\
&=&\left(\bigoplus_{E_i^M\Gamma_M\in\Gamma/\Gamma_M}E_i\right)M^{n_j}\\
&=&\bigoplus_{E_i^M\Gamma_M\in\Gamma/\Gamma_M}E_iM^{n_j}
\end{eqnarray*}

However, we have
$$\widehat{T[\Gamma_M]}\widehat{T[\Gamma_M^0]}\widehat{S[\Gamma_M^0]}M_j
\simeq \widehat{T[\Gamma_M]}
\bigoplus_{E_i\Gamma_M\in\Gamma/\Gamma_M}E_i^M(M_j).$$
The right hand side is a direct factor of
$\widehat{T[\Gamma_M]}\widehat{T[\Gamma_M^0]}(X_j)$.
We get that
$\widehat{T[\Gamma_M]}\widehat{T[\Gamma_M^0]}\widehat{S[\Gamma_M^0]}M_j$ is a direct factor
of $\widehat{T[\Gamma_M]}\widehat{T[\Gamma_M^0]}(X_j)$. But $\widehat{S[\Gamma_M^0]}M_j=X_j\oplus Y_j$
and therefore
$$\widehat{T[\Gamma_M]}\widehat{T[\Gamma_M^0]}\widehat{S[\Gamma_M^0]}M_j
=\widehat{T[\Gamma_M]}\widehat{T[\Gamma_M^0]}
X_j\oplus \widehat{T[\Gamma_M]}\widehat{T[\Gamma_M^0]}Y_j.$$
By the Krull-Schmidt property of $\mathcal H$, we conclude
$\widehat{T[\Gamma_M]}\widehat{T[\Gamma_M^0]}(Y_j)=0$.
But $T=T[\Gamma_M]T[\Gamma_M^0]=\widehat{T[\Gamma_M]}\widehat{T[\Gamma_M^0]}D$ by
Lemma~\ref{monadonorbitcategory}. Now,
$Y_j$ is a direct factor of some $D(\widetilde Y_j)$. More precisely,
$Y_j=(\widetilde Y_j,e)$ for some non trivial idempotent endomorphism $e$ of
$\widetilde Y_j$.
Since by Proposition~\ref{ehatmapstoisomorphiccopies} the counit $ST\lra \textup{id}_{\mathcal G}$ is
pointwise split, $\widetilde Y_j$ is a direct
factor of $ST\widetilde Y_j$. Since $\widehat{T[\Gamma_M]}\widehat{T[\Gamma_M^0]}(Y_j)=0$
we get $Te=\widehat{T[\Gamma_M]}\widehat{T[\Gamma_M^0]}D(e)=0$ which is a contradiction.
This proves the theorem.
\dickebox

\begin{Rem} \label{classicalcase}
Suppose that $G$ is a group with finite index normal subgroup $N$.
Then let ${\mathcal H}=kH-mod$
for a commutative ring $k$, and ${\mathcal G}=kG-mod$, the functors $S$ and $T$
being induction and restriction. The counit of the adjunction
$ST\lra \textup{id}_{\mathcal G}$ is then the trace map, as was developed by Linckelmann~\cite{Linckelmanntransfer}.
Proposition~\ref{monadsadmittingsection} requires this map to be split, which is
equivalent to $|\Gamma|=|G:H|$ being invertible in $k$.  In this case all
adjunctions giving this monad are then actually isomorphic.
Theorem~\ref{CliffordforMonads} hence generalises
Clifford's theorem~\ref{Cliffordtheoremclassic} for group representations for
normal subgroups with index invertible in the ground ring.
\end{Rem}

\begin{Rem}
Recently Asashiba developed his paper \cite{AsashibaGalois1} further. He replaced
the action of the group $\Gamma$ by the Grothendieck
construction~\cite{AsashibaGrothendieck}. We are not
yet able to generalise Theorem~\ref{CliffordforMonads}
in this direction but we intend to do so in future work.
\end{Rem}

\section{Examples and applications}

\label{Exampleapplicationsection}

\subsection{Galois modules}

In Remark~\ref{classicalcase} we observed that the classical Clifford theorem over fields of
characteristic $0$ is a special case of
Theorem~\ref{CliffordforMonads}.

We consider a slight generalisation using Galois modules.

Let $K$ be a field and suppose that $K\leq L\leq M$ is a sequence of finite field extensions.
\begin{itemize}
\item
Suppose that $M$ is Galois over $K$, denote $Gal(M:K)=:\Gamma$ and
\item suppose that $L$ is normal, i.e.
$\Delta:=Gal(M:L)\unlhd Gal(M:K)=\Gamma$. Then, restriction gives an epimorphism
$$Gal(M:K)\lra Gal(L:K)$$ with kernel $\Delta=Gal(M:L)$.
\item Suppose that $G$ is a finite group and
let $\varphi:G\lra\Gamma$ be a group homomorphism. Then, via $\varphi$ we get that
$M$ is a $KG$-module, and via $$\Gamma\lra \Gamma/\Delta\simeq Gal(L:K)$$
we obtain that also $L$ is a $KG$-module.
\item
Let $H\unlhd G$ be a normal
subgroup, and suppose further that $[\varphi(H),\Delta]=1$.
\end{itemize}
Consider the twisted group rings $M\sdp G$, respectively $L\sdp G$.
Recall that $M\sdp G$ is $M\otimes_KKG$ as abelian group, and the multiplicative law is given by
$$(m_1\otimes g_1)\cdot(m_2\otimes g_2)=(m_1\cdot \varphi(g_1)(m_2))\otimes g_1g_2$$
for any $m_1,m_2\in M$ and $g_1,g_2\in G$.

Then $L\sdp H$ is a subring of $M\sdp G$.
This situation then allows to apply our Theorem~\ref{CliffordforMonads}.

\begin{Lemma}\label{twistedgroupringsrestriction}
With the notations above, $M\sdp G$ is a free $L\sdp H$-module with basis parameterised by
$\Delta\times G/H$.
\end{Lemma}

Proof. Since $M$ is Galois over $L$, by the normal basis theorem we get
$$M=\bigoplus_{\delta\in\Delta}\ ^\delta L.$$
Hence, since $KG=\bigoplus_{Hg\in G/H}KHg$,
\begin{eqnarray*}M\sdp G&\simeq&\left(\bigoplus_{\delta\in\Delta}\ ^\delta L\right)\sdp G\\
&\simeq&\bigoplus_{\delta\in\Delta}\left(\ ^\delta L\sdp G\right)\textup{ since $L$ is normal
and hence $G$ acts on $L$}\\
&\simeq&\bigoplus_{\delta\in\Delta}\bigoplus_{Hg\in G/H}\left(\ ^\delta L\sdp H\right)g\\
&\simeq&\bigoplus_{(\delta,Hg)\in\Delta\times G/H}\left(\ ^\delta L\sdp H\right)g
\end{eqnarray*}
Since $\delta\in\Delta,$ the automorphism $\delta$ fixes $L$ pointwise, and hence
$$M\sdp G\simeq \left(L\sdp H\right)^{|\Delta\times G/H|}$$
as $L\sdp H$-module.  \dickebox

\medskip

$G$ acts on $\Gamma$ via $\varphi$, and then by conjugation. We may hence form the semidirect product
$\Gamma\sdp G$. As $L$ is normal over $K$, we get $\Delta\unlhd \Gamma$, and therefore by
restriction we get
a semidirect product $\Delta\sdp G$. Further, $[\varphi(H),\Delta]=1$ implies that this induces a
semidirect product structure $\Delta\sdp G/H$.

\begin{Lemma}  \label{normalsubgroupsituationGaloismodule}
The monad given by the adjoint pair $(M\sdp G\otimes_{L\sdp H}-,res^{M\sdp G}_{L\sdp H})$
is a normal subgroup situation with group $\Delta\sdp G/H$.
\end{Lemma}

Proof. Denote $(\Psi,\Phi):=(M\sdp G\otimes_{L\sdp H}-,res^{M\sdp G}_{L\sdp H})$  for short.
From Lemma~\ref{twistedgroupringsrestriction} we see that
$$\Phi\circ\Psi=\bigoplus_{(\delta,g)\in \Delta\sdp G/H}(\ ^\delta L\sdp H)g.$$
Hence we may put
$$E_{(\delta,g)}=(\ ^\delta L\sdp H)g\otimes_{L\sdp H}-.$$
Composition gives
$$E_{(\delta_1,g_1)}\circ E_{(\delta_2,g_2)}=E_{(\delta_3,g_3)}$$
for $(\delta_3,g_3)=(\delta_1\cdot\ ^{\varphi(g_1)}\delta_2\;,\; g_1g_2)$.
This is precisely the structure of the group $\Delta\sdp G/H$.
\dickebox

\begin{Rem}
Note that we did not fully use that $K, L$ and $M$ are actually fields.
As long as we have a normal basis theorem, such as for tame abelian extensions of algebraic
integers in number fields (using Taylor's theorem~\cite{Taylor}), the statement stays true.
\end{Rem}

We now consider the functor $\Psi:=(M\sdp G)\otimes_{(L\sdp H)}-$
$$L\sdp H-mod\lra M\sdp G-mod$$
together with its right adjoint, the forgetful functor $\Phi$.
By Lemma~\ref{normalsubgroupsituationGaloismodule}
the adjoint pair $(\Psi,\Phi)$ determines the monad
$$\Phi\circ \Psi=\bigoplus_{(\gamma,gH)\in\Delta\sdp G/H}\ E_{(\delta,g)}$$
where $$E_{(\delta,g)}=(\ ^\delta L\sdp H)g\otimes_{L\sdp H}-.$$

Hence, the hypotheses of Theorem~\ref{CliffordforMonads}
are satisfied.

A particularly interesting case occurs when
$K$ is a field of characteristic $0$. Indeed, again the
counit $$\Psi\circ\Phi\lra\textup{id}$$
is the trace function. The hypothesis that $K$ is a field of characteristic $0$ then
shows that the trace function, hence the counit of the adjunction, is split.
Therefore, by Proposition~\ref{monadsadmittingsection} all adjunctions are isomorphic.
Hence, the theorem gives a criterion how $M\sdp G\otimes_{L\sdp H}V$
decomposes in terms of the inertia group $I_{\Delta\sdp G/H}(V)$
for a simple $L\sdp H$-module $V$. This inertia group also depends on
the Galois extension. However, $\Delta=Gal(M:L)$, and hence $^\delta L\simeq L$ as $L$-module
for any $\delta\in\Delta$. Therefore $I_{\Delta\sdp G/H}(V)=\Delta\sdp U$ for some subgroup $U$ of $G$
containing $H$.

\subsection{De-equivariantization}

In the theory of tensor, braided and fusion categories the concept studied in our paper is known
under the name of de-equivariantization. Recall the setting from \cite[Section 8.22]{EtingofGelaki}
there. For more ample details we refer to \cite{EtingofGelaki}.

A fusion category is a finite semisimple
tensor category such that the endomorphism algebra of the unit element is the base field $K$.
Let ${\mathcal C}$ be a fusion category with an action of a group $G$. Then $KG-mod$ is a
subcategory of the full subcategory of $G$-equivariant objects in ${\mathcal C},$ which is again a
fusion category denoted ${\mathcal C}^G$. This embedding factors through the canonical embedding
$Z({\mathcal C}^G)\hookrightarrow{\mathcal C}^G$.

Let $\mathcal D$ be a fusion category and let $G$ be a finite group acting on $\mathcal D$.
Suppose that there is a braided tensor functor $KG-mod\lra Z({\mathcal D})$ such that the composition
with $ Z({\mathcal D})\hookrightarrow {\mathcal D}$ is fully faithful. Then
the image of the regular $KG$-module in ${\mathcal D}^G$ forms a monad $A$. Consider
$A-Mod_{\mathcal D}$ (cf Section~\ref{EilenbergMooreSection}). Then ${\mathcal D}_G:=A-Mod_{\mathcal D}$
is again a fusion category with an action of $G$ such that $({\mathcal D}_G)^G\simeq{\mathcal D}.$
Observe that the monad $A$ is precisely what we call a normal subgroup situation. We hence may
apply the results of Theorem~\ref{CliffordforMonads} to this setting.

\bigskip

\subsection*{Acknowledgement:} I wish to thank Javad Assadollahi for having invited me to the
CIMPA UI IPM Isfahan school and conference on representations of algebras in April 2019  in
Isfahan, Iran and to Hideto Asashiba with whom I had the chance to discuss topics
on his lecture~\cite{AsashibaGrothendieck}. This lecture
was the origin of the present paper. I thank Benjamin Sambale for having asked how to
interpret Galois theory in terms of Clifford theory, and to the anonymous referee for
having asked for supplementary applications, which lead to Section~\ref{Exampleapplicationsection}.
I thank Alexis Virelizier for mentioning the results from \cite[Section 8.22]{EtingofGelaki}.

\end{document}